\documentclass[12pt]{amsart}

\newcommand{\Z}{{\mathbb Z}}
\newcommand{\R}{{\mathbb R}}
\newcommand{\C}{{\mathbb C}}

\begin{document}

\title{COMPLEXES OF CONNECTED GRAPHS}

\author{V.A. Vassiliev}

\email{vva@mi.ras.ru}

\date{Revised version published in 1993}

\begin{abstract}
Graphs with the given k vertices generate an (acyclic) simplicial complex. We
describe the homology of its quotient complex, formed by all connected graphs,
and demonstrate its applications to the topology of braid groups, knot theory,
combinatorics, and singularity theory. The multidimensional analogues of this
complex are indicated, which arise naturally in the homotopy topology, higher
dimensional Chern-Simons theory and complexity theory.\end{abstract}

\maketitle

\section{Introduction}

Let $A$ be a set of $k$ elements. Denote by $\Delta(A)$ the simplex with
$\binom{k}{2}$ vertices, being in one-to-one correspondence with the
two-element subsets of $A$. Any face of this simplex can be depicted by a graph
with $k$ vertices corresponding to the elements of $A$: this graph contains the
segment between some two elements $a, b \in A$ iff the vertex corresponding to
the set $(a, b)$ belongs to the chosen face.

A graph is called {\em connected} if any two points in $A$ can be joined by a
chain of its segments.

Denote by $M(A)$ the set of all faces in $\Delta(A)$ corresponding to all
non-connected graphs. Obviously, this is a subcomplex of the simplex
$\Delta(A)$.

\medskip
\noindent  {\bf Definition.} The {\it complex of connected graphs} associated
with $A$ is the quotient complex $\Delta(A)/M(A)$; the notation of this
quotient complex is $K(A)$.
\medskip

Of course, all these objects corresponding to different $k$-element sets $A$
are isomorphic. We shall denote by $\Delta_k$, $M_k$ and $K_k$ the complexes
$\Delta(A)$, $M(A)$ and $K(A)$ where $A$ is the set of naturals $1, 2, ... ,
k$.

\medskip \noindent
{\bf Theorem 1}. {\it For any $k$-element set $A$, and coefficient group $G$,
the homology group $H_i(K(A),G)$ is trivial for $i \neq k - 2$, and}
\begin{equation}\label{one1} H_{k-2}(K(A), G) = G^{(k-1)!} . \end{equation}

\medskip
This theorem is a special case of a general theorem of Folkman [F]. We give an
independent proof, based on the geometrical considerations, namely, on the
theory of plane arrangements.

Here is an explicit realization of the group $H_{k-2}(K(A))$. Let us
distinguish one point $* \in  A$.

\medskip \noindent
{\bf Theorem 2}. {\it For a basis in the group $H_{k-2}(K (A))$ we can take the
classes of $(k - 1)!$ trees, homeomorphic to a segment $($i.e., having only
vertices of orders 1 or 2$)$, and such that an endpoint of this segment $($a
vertex of order 1$)$ coincides with the distinguished point} $* \in A$.
\medskip

Denote by $\nabla(A)$ the simplex with $2^{k-1} - 1$ vertices, being in
one-to-one correspondence with the nonordered partitions of $A$ into 2 nonempty
subsets. A face of this simplex is called {\em non-complete} if there exist two
elements $a, b \in A$ such that, for any vertex of this face, $a$ and $b$ lie
in the same part of the corresponding partition of $A$. The set of all
non-complete faces is, obviously, a subcomplex in $\nabla(A)$; we shall denote
this complex by $\Omega(A)$.

\medskip \noindent
{\bf Definition}. The {\it complex $C(A)$ of complete partitions} associated
with $A$ is the quotient complex $\nabla(A)/\Omega(A)$.
\medskip

Again, the complexes $\nabla_k$, $\Omega_k$, $C_k$ are the complexes
$\nabla(A),$ $\Omega(A)$, $C(A)$ for $A = (1,2, ..., k)$.

\medskip \noindent
{\bf Theorem 3}. {\it For any finite set A and any coefficient group G, there
is a natural isomorphism}
\begin{equation} \label{two2} H_*(C (A), G) = H_*(K(A), G).\end{equation}

\medskip
The number $(k-1)!$ from the formula (1) is well-known to the specialists in
the braid groups. Indeed, the Poincar\'e polynomial of the cohomology of the
colored braid group with $k$ strings, $I(k)$, equals
\begin{equation}\label{three3}(1 + t)(1 + 2t) \cdots (1 + (k-1)t);\end{equation}
in particular $\mbox{dim }H_{k-1}(I(k)) = (k-1)!$, see [A1]. This is not an
occasional coincidence: the complex of connected graphs appears naturally in
the description of the homotopy type (and, in particular, of the homology
groups) of the classifying space of the group $I(k)$. This description depends
on the fact, that the classifying space $K(I(k), 1)$ can be represented as the
complement to a certain collection of planes in the space ${\mathbb C}^k$, and
on a general formula (see [GM], EV11, [V2], [ZZ]) expressing the stable
homotopy type of the complement of arbitrary collection of affine planes in
${\mathbb R}^n$ in terms of the combinatorial characteristics of this
collection (namely, of the dimensions of the planes and all their
intersections).

Indeed, let $V_1,   \dots , V_s$ be a collection of affine planes in ${\mathbb
R}^n$ (such collections are called {\it affine plane arrangements}), and $V$
the union of these planes. Let $\Lambda$ be any affine plane which is the
intersection of several planes $V_i$. Consider the simplex $\Delta(\Lambda)$
whose vertices correspond formally to all the planes $V_i$ containing
$\Lambda$. A face of this simplex, i.e., a collection of such planes, is called
{\it marginal}, if the intersection of these planes is strictly greater than
$\Lambda$. The quotient space of the simplex $\Delta(\Lambda)$ by the union of
all marginal faces will be denoted by $K(\Lambda)$.

\medskip \noindent
{\bf Theorem 4} (see [V1], [V2], [ZZ]). {\it The one-point compactification of
the variety $V$ is homotopy equivalent to the wedge of $\mbox{\rm
dim}(\Lambda)$-fold suspensions of the spaces $K(\Lambda)$ taken over all
planes $\Lambda$ which are the intersections of some planes of the set $(V_1,
\dots, V_s)$}.

\medskip \noindent
{\bf Corollary 1}. {\it The stable homotopy type of the complement ${\mathbb
R}^n - V$ of the arrangement $V$ is completely determined by its combinatorial
$($dimensional$)$ characteristics}.
\medskip

Indeed, this complement is Spanier-Whitehead dual to the one-point
compactification of $V$; hence the stable homotopy types of these spaces are
completely determined one through the other; see [W].

\medskip \noindent
{\bf Corollary 2}. {\it The cohomology group $H^i({\mathbb R}^n - V)$ is
isomorphic to

\begin{equation} \label{four4} \bigoplus \tilde
H_{n-i-1-\mbox{\rm \footnotesize dim}(\Lambda)}(K(\Lambda)\end{equation} $($the
summation over all planes $\Lambda$ which are the intersections of several
planes $V_i)$.}
\medskip

Indeed, this follows from Theorem 4 by the Alexander duality.

The formula (4) was previously obtained by Goresky and MacPherson, see [GM].

\medskip \noindent
{\bf Corollary 3}. {\it For any $k$-element set $A$, the homology group
$H_i(K(A),G)$ is trivial for $i \neq k - 2$; see Theorem 1.} \medskip

Indeed, the lower estimate follows from the fact that any connected graph with
$k$ vertices has at least $k-1$ segments. To prove the upper estimate, consider
the hyperplane arrangement in ${\mathbb R}^k$, given by all ``diagonal'' planes
distinguished by the equalities $x_i = x_j$, $1 \leq i < j \leq k$. Let
$\Lambda$ be the main diagonal of this arrangement, i.e., the line $\{x_1 =
\dots = x_k\}$. The complex $K(\Lambda)$ is naturally isomorphic to the complex
$K_k$. On the other hand, by Theorem 4, the homology group $H(K(\Lambda))$ is a
direct summand in the $(l+1)$-dimensional homology group of the one-point
compactification of our arrangement in ${\mathbb R}^k$. Since this
compactification is a $(k-1)$-dimensional complex, the desired estimate
follows.

Another important application of the complex $K(A)$ appears in knot theory. In
[V2]--[V4], I have constructed a spectral sequence providing the cohomology
classes of the spaces of knots (i.e., of nonsingular embeddings $S^1 \to
{\mathbb R}^n$, $n \ge  3$). The complex $K(A)$ is an essential part of this
spectral sequence; see Section 4.

The link theory provides a similar complex: {\it the complex of connected
graphs with colored vertices}.

Suppose that the set $A$ of the vertices of our graphs is partitioned into $d$
parts: $A = A_l \cup \dots \cup A_d$, $\mbox{card} A_i = k_i$, $\rm{card} A =
k_l + \dots + k_d = k$. A graph with the vertices in A is called {\em
concordant with the partition} if the endpoints of any of its segments belong
to different parts of this partition.

Again, all the concordant graphs generate an (acyclic) complex which is
naturally isomorphic to the simplicial complex of the obvious triangulation of
a simplex with $k_1 k_2 + k_1 k_3 + \dots + k_{d-1}k_d$ vertices. The connected
concordant graphs (i.e., the concordant graphs connecting all points of the set
$A$) constitute a quotient complex of this complex. We denote this quotient
complex by $K(A; A_1, \dots, A_d)$.

The usual complex $K(A)$ is a special case of this one, which corresponds to
the partition into separate points.

\medskip \noindent
{\bf Theorem 5}. {\it The complex $K(A; A_1,... , A_d)$ is acyclic in the
dimensions not equal to $k - 2$.}
\medskip

The proof is almost the same as that of Corollary 3. For more about the
topological applications of this complex, see Section 6.
\medskip

There exist natural multidimensional analogues of the complexes of connected
graphs: the complexes $K(A, t)$, $t > 2$, playing the same role in the higher
dimensional Chern-Simons theory and equations of multidimensional simplices
(see [FNRS], [MS]) which the standard complexes $K(A)$ and $K(A; A_1, \dots,
A_d)$ play (by means of the braid and knot theories) in the usual variants of
these theories. These complexes were for the first time investigated by A.
Bj\"orner and V. Welker [BW] in a connection with the problems of the
complexity theory; see also [BLY].

The complex $K(A)$ considered above appears in [V2]--[V4] from the resolution
of the discriminant variety, i.e., of the (closure of) set of all maps $S^1 \to
{\mathbb R}^3$ having self-intersections. The ``self-intersection'' is a
bisingularity: it is defined by a condition imposed on a map at two different
points of the issue manifold. Now suppose that we define the discriminant by a
condition imposed at $r$ points, $r \geq 3$ (the simplest problem where such
discriminants arise is the study of the space of maps $S^1 \to {\bf R}^n,$ $n
\geq 2$, having no triple self-intersections, or, more generally, the
self-intersections of multiplicity $r$). The homology groups of such spaces are
provided by a spectral sequence similar to the one from [V2]--[V4]; but the
role of the complex $K(A)$ here is played by the complexes $K(A, t)$ described
as follows.

For any $t \geq 2$ consider the simplex $\Delta(A, t)$ with $\binom{k}{t}$
vertices being in one-to-one correspondence with the subsets of cardinality $t$
in $A$. A face of this simplex, i.e., a collection of such subsets, is called
{\it connected}, if any two points in $A$ can be joined by a chain in $A$, any
two neighboring elements of which belong to some subset of our collection.

Again, the connected faces form a quotient complex in $\Delta(A, t)$; denote it
by $K(A, t)$. For $t = 2$, this is exactly the usual complex of connected
graphs considered above.

\medskip \noindent
{\bf  Theorem 6} (see [BW]). {\it Suppose that $t > 2$. Then the homology
groups $H_i(K(A,t))$ are trivial for all $i$ except maybe for $i$ of the form
$k-1-q(t-2)$, $q \geq 1$.}
\medskip

(The fact that $H_i(K(A(t)) = 0$ for $i > k - (t - 1)$ can be proved again in
the same way as Corollary 3 to Theorem 4.)

The complex $K(A, t)$ has also the versions with ``colored vertices'': the
simplest (and very absorbing) problem where such complexes appear is the
topological classification of triplets of closed curves in ${\mathbb R}^2$
having no common points.

I am grateful to V.I.~Arnold, A.~Bj\"orner, I.M.~Gelfand, M.M.~Kapranov,
M.L.~Kontsevich, G.L.~Rybnikov, B.Z.~Shapiro, V.V.~Schechtman and G.~Ziegler
for helpful discussions.

\section{Complex of connected graphs and complex of complete partitions. Proofs
of Theorems 1, 2, 3}

\noindent {\bf Proof of Theorem 3.}

\medskip \noindent
{\bf Lemma 1}. {\it For any $i \geq 0$ and any coefficient group $G,$ $\tilde
H_i(M(A),G) = H_{i+1}(K(A),G),$ $\tilde H_i(\Omega(A),G) = H_{i+1}(C(A),G)$,
where $\tilde H$ denotes the homology reduced modulo a point.}
\medskip

This follows immediately from the acyclicity of simplices.

Therefore we have only to prove that $H_*(M(A),G) = H_*(\Omega(A),G)$. To do
this, we construct a cellular complex
$$\sigma(A) \subset M(A) \times \Omega(A)$$
which is homology equivalent to both $M(A)$ and $\Omega(A)$ (here we do not
distinguish between the complexes $M(A)$, $\Omega(A)$ and the topological
spaces defined as the unions of corresponding faces of the simplices).

Consider an arbitrary partition $\pi$ of $A$ into $n$ nonempty subsets, $2 \leq
n \leq k -1$. Let $\alpha(\pi)$ be the simplex in $M(A)$ whose vertices
correspond to all two-element subsets $(a, b) \subset A$ such that $a$ and $b$
are in the same part of the partition. Let $\beta(\pi)$ be the simplex in
$\Omega(A)$, whose vertices correspond to all partitions of $A$ into two
nonempty subsets such that it is a subpartition of these partitions. The
product $\alpha(\pi) \times \beta(\pi)$ is a cell in $M(A) \times \Omega(A)$.
Define the subset $\sigma(A) \subset M(A) \times \Omega(A)$ as the union of
such products $\alpha(\pi) \times \beta(\pi)$ over all nontrivial partitions
$\pi$ of the set $A$. Now, Theorem 3 follows from the following lemma.
\medskip

\noindent {\bf Lemma 2.} {\it The obvious projections of the set $\sigma(A)$
onto $M(A)$ and $\Omega(A)$ induce the isomorphisms of the homology groups}.
\medskip

\noindent {\bf Proof.} Consider the increasing filtration of the complex $M(A)$
by its skeletons, and the filtration of $\sigma(A)$ by the preimages of these
skeletons under the obvious projection $\sigma(A) \to  M(A)$. Over any open
simplex in $M(A)$, this projection is a trivial fiber bundle whose fiber is a
simplex. Therefore, the homological spectral sequences, generated by these two
filtrations, are isomorphic beginning with the term $E^l$ and this isomorphism
is induced by the projection. Thus, $H_*(\sigma(A)) = H_*(M(A))$. For the
complex $\Omega(A)$ the proof is exactly the same, and Lemma 2 is proved.
\medskip

\noindent {\it Remark}. In fact the projections of $\sigma(A)$ onto $M(A)$ and
$\Omega(A)$ induce also the homotopy equivalences of all these spaces.
\medskip

\noindent {\bf Proof of Theorem 1}. By Corollary 3 and Theorem 3, this theorem
is equivalent to the following:

\noindent {\bf Lemma 3}. {\it The Euler characteristic of the complex
$\Omega(A)$ is equal to} $1- (-1)^k (k - 1)!$.
\medskip

\noindent {\bf Proof of Lemma 3}. Let $F_p$ be the union of all simplices
$\beta(\pi) \subset \Omega(A)$ considered in the proof of Theorem 3 and such
that $\pi$ is a partition of $A$ into $n$ parts, $n \leq p$. The sets $F_p$
define an increasing filtration of the space $\Omega(A)$, $F_2 \subset F_3
\subset \dots \subset F_{k-1} = \Omega(A)$. Consider the spectral sequence
$E^r_{p,q}$ calculating the homology of $\Omega(A)$ and generated by this
filtration.

By definition, $E^1_{p,q} = H_{p+q}(F_p, F_{p-1}).$
\medskip

\noindent {\bf Sublemma}. {\it For any $p = 2,3,..., k -1$, the group
$E^1_{p,q}$ splits into a direct sum of subgroups corresponding to different
partitions of $A$ into $p$ nonempty subsets, and any of these groups is
isomorphic to the group} $H_{p+q}(C_p) = H_{p+q-1}(\Omega_p ; \mbox{a point})$.
\medskip

Consider any such partition $\pi$ of $A$ and the corresponding simplex
$\beta(\pi) \subset \Omega(A)$. This simplex has $2^{p-1} - 1$ vertices
corresponding to all partitions of $A$ into two subsets such that $\pi$ is
subordinate to these partitions; thus this simplex can be identified with the
simplex $\nabla([\pi])$ where $[\pi]$ is the set of parts of the partition
$\pi$. A face of our simplex $\beta(\pi)$ does not belong to $F_{p-1}$ iff the
corresponding face in $\nabla([\pi])$ corresponds to a complete partition; in
this case this face cannot be a face of some other simplex $\beta(\pi')$ where
$\pi'$ is a partition of $A$ into $p$ parts, $\pi' \neq \pi$, and the sublemma
follows. \medskip

Now suppose that Lemma 3 is proved for all complexes $\Omega(B)$ with
$\mbox{card}(B) < k = \mbox{card}(A)$. Then it follows from the sublemma that
the wanted Euler characteristic equals
$$ \sum_{n=2}^{k-1} (-1)^n (n-1)! \langle k ? n \rangle , $$
where $\langle k?n\rangle$ is the number of nonordered partitions of a
$k$-element set into $n$ nonempty parts. Since $\langle k?1\rangle = \langle
k?k\rangle = 1$, Lemma 3 is equivalent to the equality
$$ \sum_{n=1}^{k} (-1)^n (n-1)! \langle k ? n \rangle =0. $$

This equality follows immediately from the obvious combinatorial identity
$$ \langle k?n\rangle = n \langle k - 1?n \rangle + \langle k - 1 ? n-1\rangle,$$
and Theorem 1 is completely proved.
\medskip

\noindent {\bf Proof of Theorem 2}. We shall prove this theorem for the complex
$K_k$, so that $A = (1, 2, ... , k)$.
\medskip

\noindent {\bf Definition.} A $k$-tree is a tree with $k$ vertices $(1), \dots,
(k)$. A k-tree is called {\it normed} if its vertex (1) is of order 1, and is
called {\it linear}, if it is homeomorphic to a segment (i.e., the orders of
all its vertices are no more than 2).
\medskip

\noindent {\bf Lemma 4.} {\it Any $k$-tree is homologous in the complex $K_k$
to a linear combination of normed $k$-trees.}

\medskip \noindent
{\bf Lemma 5}. {\it Any normed $k$-tree is homologous in the complex $K_k$ to a
linear combination of normed linear trees.}
\medskip

These two lemmas together with Theorem 1 imply Theorem 2, because the number of
all linear normed k-trees equals exactly $(k - 1)!.$
\medskip

\noindent {\bf Proof of Lemma 4}. Consider any two segments of our tree having
the common vertex (1). Consider a graph with $k$ segments obtained from our
tree by adding the segment joining the other endpoints of these two segments.
Obviously, the boundary of this graph (more rigorously, of the face,
corresponding to this graph) is a linear combination of our tree, of two trees
having one less order of the vertex (1), and $k-3$ graphs containing a triangle
(and hence non-connected). This implies the lemma.
\medskip

\noindent {\bf Proof of Lemma 5}. Let us start from the vertex (1) of our
normed tree and go along this tree until we meet for the first time a vertex of
order $> 2$. (Before this moment our path is uniquely determined). Choose any
two segments of this graph having the endpoints at this vertex and not
coinciding with the segment by which we came to it. Consider the graph with $k$
segments obtained from our tree by adding the segment connecting two other
endpoints of these two segments. Again, the boundary of this graph is
homologous in the complex $K_k$ to the linear combination of our $k$-tree and
two trees with one less order of our vertex. This proves Lemma 5.

\section{The topology of the complement of a plane arrangement.
Proof of Theorem 4. Colored braid group}

\noindent {\bf 3.1. Notations.} Let $V_1,\dots , V_s$ be a finite set of affine
planes in $\R^n$, and $V$ the union of all $V_i$. Denote by $S$ the set of
naturals $1, \dots, s$. For any subset $J \subset S$, $V_J$ is the intersection
of planes $V_j$, $j \in J$. The dimension of this intersection is denoted by
$|J|$. For any $J$, $J'$ is the maximal subset in $S$ such that $V_J = V_{J'}$.
$\bar V_J$ and $\bar V$ are the notations for the one point compactifications
of $V_J$ and $V$.

A set $J$ is called {\it geometrical} if $J = J'$.

For any geometrical set $J \subset S$, denote by $\Delta(J)$ the simplex whose
vertices are in one-to-one correspondence with the elements of $J$. (In the
Introduction this simplex was denoted by $\Delta(V_J)$: this change cannot lead
to a misunderstanding.) Let $M(J)$ be the subcomplex in $\Delta(J)$ consisting
of all marginal faces, and recall the notation $K(V_J)$ for the quotient
complex $\Delta(J)/M(J)$.
\medskip

\noindent {\bf 3.2. Geometrical resolution of the complex $V$.} Without loss of
generality, we shall assume that the dimensions of all planes $V_j$ are
positive. Consider a space $\R^N$, where $N$ is sufficiently large, and some
$s$ affine embeddings $I_j: V_j \to \R^N$, $j \in S$. For any point $x \in V$,
consider all its images in $\R^N$ under all maps $I_j$ such that $x \in V_j$.
Denote by $\# x$ the number of such $j$ and by $x'$ the convex hull of these
images in $\R^N$.
\medskip

\noindent {\bf Lemma 6.} {\it If $N$ is sufficiently large and the system of
embeddings $I_j$ is generic, then for any point $x \in V$ the polyhedron $x'$
is a simplex with $\# x$ vertices, and the intersection of simplices $x'$, $y'$
is empty if $x \neq y$.}
\medskip

The proof is trivial.
\medskip

We shall suppose that the maps $I_j$ satisfy this lemma.

Denote by $V'$ the union of all simplices $x'$ over all $x \in V$ and by $\bar
V'$ the one-point compactification of $V'$. These spaces $V', \bar V'$ will be
called the {\it geometrical resolutions} of $V$ and $\bar V$. The natural
projection $\pi : V' \to V$ (which maps any simplex $x'$ into the point $x$) is
obviously proper and can be extended to a continuous map $\bar V' \to \bar V$
which will be denoted by the same letter $\pi$.
\medskip

\noindent {\bf Lemma 7}. {\it The projection $\pi: \bar V' \to \bar V$ induces
a homotopy equivalence of these spaces.}
\medskip

This is (a special case of) the principal fact of the theory of simplicial
resolutions, see f.i. [D].

Hence we have only to prove the following theorem.
\medskip

\noindent {\bf Theorem $4'$}. {\it For any affine plane arrangement $V$, the
one point compactification of its resolution $V'$ is homotopy equivalent to the
wedge indicated in Theorem 4.}
\medskip

\noindent {\bf 3.3. Proof of Theorem $4'$} This proof is based on a variation
of stratified Morse theory, see [GM].
\medskip

\noindent {\bf Definition}. A function $f : \R^n \to \R^1$ is called a {\it
generic quadratic function} if it can be expressed in the form $x_1^2 + ... +
x_n^2$ in some affine coordinate system in $\R^n$, whose origin does not belong
to $V$, and any level set $f^{-1}(t)$ of this function is tangent to at most
one plane $V_J$. \medskip

Let us fix such a function $f$. For any set $J$, denote by $t_J$ the only
number $t$ such that $f^{-1}(t)$ is tangent to $V_J$; all numbers $t_J$ are
called {\em singular values}, and the other values are {\em regular}.

For any value $t \in \R^1$ denote by $V'(t)$ the space

$$\pi^{-1} (V \cap f^{-1}([t, \infty)).$$

\medskip \noindent
{\bf Lemma 8}. {\it $(a)$ If $t$ is greater than all singular values $t_J$,
then the quotient space $V'/ V'(t)$ is homotopy equivalent to $\bar V'$;
\smallskip

$(b)$ If $t$ is less than all values $t_J$, then this quotient space is a
point;
\smallskip

$(c)$ if the segment $[t, s]$ does not contain singular values, then the
obvious factorization mapping $(V'/V'(s)) \to (V/V'(t))$ is a homotopy
equivalence;
\smallskip

$(d)$ if the segment $[t, s]$ contains exactly one singular value $t_J$, $t <
t_J < s$, and the set $J$ is geometrical, then the space $V'/V'(s)$ is homotopy
equivalent to the wedge of spaces $V'/ V'(t)$ and
$\Sigma^{|J|}(\Delta(J)/M(J))$, where $\Sigma^i$ is the notation of the
$i$-fold reduced suspension.}
\medskip

Theorem $4'$ follows immediately from this lemma.

Items a, b and c of this lemma are obvious; let us prove d. For any geometrical
set $J$, define the {\it proper inverse image} of the plane $V_J$ as the
closure in $\R^N$ of the union of simplices $x'$ over all points $x \in V_J$
which do not belong to subplanes $V_I$ of lower dimensions in $V_J$. Denote
this closure by $V'_J$, and by $V'_J(t)$ its intersection with $V'(t)$.

Any space $V'_j$ is naturally homeomorphic to the complex $\Delta(J) \times V_J
\simeq \Delta(J) \times \R^{|J|}$.

Now let $J$ be the geometrical set considered in Lemma 8(d). Then, the space
$V'_J/V'_J(s)$ is naturally homeomorphic to the space $$(\Delta(J) \times V_J)/
(\Delta(J) \times (V_J \cap f^{-1}([s, \infty))) \simeq
\Sigma^{|J|}\Delta(J).$$

Let $W$ be the union of preimages of all other planes $V_I$, $I \neq J$, having
nonempty intersections with the disk $f^{-1}([0, s])$. The intersection of
varieties $V'_J$ and $W$ can be naturally identified with the complex $M(J)
\times V_J \cong M(J) \times \R^{|J|}$. Let $\phi$ be the identical imbedding
of this intersection into W. Let $\psi$ be the map from the quotient space
$(M(J) \times V_J)/(M(J) \times V_J(s))$ into $W/(W \cap V'(s))$ induced by
$\phi$. The space $V'/V'(s)$ can be considered as the space $$(W/(W \cap
V'(s))) \cup_\psi (V'_J/V'_J(s)),$$ where $\cup_\psi$ is the topological
operation ``paste together by the map $\psi$''; see e.g. [FV].

But the quotient space $W/(W \cap V'(s))$ is naturally homotopy equivalent to
the space $V'/V'(t) = W/(W \cap V'(t))$: this homotopy equivalence is realized
by the obvious factorization map which contracts all points $z$ at which
$f(\pi(z)) \in [t, s]$.

The composition of the map $\psi$ and this factorization is a map into one
point. Since the operation $\cup_\psi$ is homotopy invariant (see [FV], section
1.2.11), the space $V'/V'(s)$ is homotopy equivalent to the composite space
$(V'/V'(t))\cup_{\psi'} (V'_J/V'_J(s))$ where the map $\psi'$ is defined on the
same subspace as $\psi$ and takes this subspace into one point $\{V'(t)\} \in
V'/V'(t)$. Hence the space $V'/V'(s)$ is homotopy equivalent to the wedge of
spaces $V'/V'(t)$ and
$$\begin{array}{l}
V'_J/(V'_J(s)\cup (M(J) \times V_J)) =  \\
= (\Delta(J) \times V_J)/((\Delta(J) \times V_J(s)) \cup (M(J) \times V_J)) = \\
= (\Delta(J)/M(J)) \wedge (V_J/V_J(s)) = \Sigma^{|J|}(\Delta(J)/M(J)).
\end{array}
$$ Q.E.D.

Theorems $4'$ and 4 are completely proved.
\medskip

\noindent {\bf 3.4. Important example: the colored braid group}
\medskip

\noindent {\bf Definition}. The {\it ordered configuration space} $F(\C^1, k)$
is the space of all ordered subsets of cardinality $k$ in $\C^1$.
\medskip

This space can be considered as a subset in $\C^k$: namely, as the complement
of the union of all planes $A_{i,j}$ distinguished by the equations $x_i =
x_j$, $i \neq j$, in the coordinates $x_1,\dots , x_k$ in $\C^k$. \medskip

\noindent {\bf Definition}. The colored braid group of $k$ strings, $I(k)$, is
the fundamental group of the space $F(\C^1, k)$. \medskip

\noindent {\bf Theorem} (see [A1]). {\it The space $F(\C^1, k)$ is a
classifying space of the group $I(k)$ : $F(\C^1, k) = K(I(k),1)$. In
particular, the cohomology of the group $I(k)$ coincides with that of the space
$F(\C^1,k)$.}
\medskip

The study of the topology of the space $F(\C^1, k)$ is a special case of the
problem considered in Theorem 4: the complex hyperplanes in $\C^k$ can be
considered as real planes of codimension 2 in $\R^{2k}$. Let us apply the
general assertion of Theorem 4 to this space. The ``deepest'' stratum of the
arrangement $\cup A_{i,j}$ is the complex line $\Lambda = (x_l = \dots = x_k)$;
it is contained in all the planes $A_{i,j}$. It is easy to see that the complex
$K(\Lambda)$ related to this stratum is exactly the complex $K_k$ of connected
graphs with $k$ vertices. Moreover, the strata of any dimension $r$ of our
arrangement are in one-to-one correspondence with the non-ordered
decompositions of the $k$-element set into exactly $r$ subsets. This implies,
in particular, the equality $H_{k-1}(I(K)) = \Z^{(k-1)!}$, as well as some
combinatorial identities which follow from the comparison of two descriptions
of the groups $H_i(I(k))$, $i < k - 1$: the first following from Theorem 4, and
the second obtained in [Al] (and expressed by formula (3)).

\section{Applications to the cohomology of spaces of knots}

In [V2]--[V4j, a system of knot invariants was constructed (and, moreover, a
way to construct the higher dimensional cohomology classes of the space of
knots was outlined). The crucial point in this construction is a simplicial
resolution of the discriminant variety, i.e., of the set of maps $S \to \R^n$
having self-intersections or points of vanishing derivative. This resolution is
a topological space together with a projection onto the discriminant variety;
the inverse image of any discriminant point $\phi$ (i.e., of a singular map
$\phi : S^1 \to \R^n$) is a simplex, whose vertices are in the one-to-one
correspondence with all the pairs of points $S^1$ glued together by the map
$\phi$, and all the points where $d\phi = 0$. In particular, over a map $\phi$
with a $k$-fold self-intersection point a simplex appears, whose vertices
correspond to the subsets of cardinality 2 in a set of $k$ points.

The space of the resolution has a natural filtration, which is defined by the
degrees of degeneracy of corresponding singular maps. For instance, the
$k$-fold self-intersection has filtration $k - 1$: coincidence of $k$ points
takes $k - 1$ independent conditions. Thus, our simplex lies in the term
$F_{k-1}$ of the filtration. A face of this simplex lies in the term $F_{k-i}$
of the filtration iff the $k$ vertex graph formed by the vertices of this face
has at least $i$ connected components, in particular it lies in $F_{k-2}$ iff
this graph is non-connected. Thus, calculating the group $\bar
H_*(F_{k-1}/F_{k-2})$ of our filtration (or, equivalently, the column
$E^1_{k-1,*}$ of the spectral sequence generated by this filtration) involves
the calculation of the homology of the complex $K_k$.

\section{Topology of the Maxwell set}

Another application of the simplicial resolutions appears from the topological
study of (the complement of ) the Maxwell set of a complex singularity.

The Maxwell set of a singularity $f : (\C^n, 0) \to (\C, 0)$ is (the closure
of) the set of parameters of a versal deformation of the singularity, which
correspond to the functions having two critical points with the same critical
value, see [AGIN]. This set can be resolved in almost the same way as the
discriminant in the knot space, and this resolution makes it possible to
calculate the homology of the complement of the Maxwell set, see [N]. Again,
the complexes $K_k$ appear naturally in these resolutions.

\section{The topological applications of the complex of the connected graphs with
colored vertices}

The complex $K(A; A_1, \dots , A_d)$ (see Theorem 5 in the Introduction)
appears naturally in the homotopy classification of links.

Indeed, consider the space $L(d)$ of all smooth maps of the disjoint union of
$d$ circles $S^1_{(1)} , \dots , S^1_{(d)}$ into $\R^3$. Obviously this set is
contractible.

Define the discriminant $\Sigma(d)$ of this space as the set of maps which send
two points of some two different circles into one point in $\R^3$.
\medskip

\noindent {\bf Definition.} A {\em link with $d$ strings} is a smooth imbedding
of the disjoint union of $d$ circles into $\R^3$. Two links are {\it homotopy
equivalent} if they lie in the same component of the complement of discriminant
in the space $L(d)$.

(The last definition is obviously equivalent to the original definition due to
Milnor [M].) In particular, the numerical homotopy link invariants are exactly
the zero-dimensional cohomology classes of the complement of the discriminant.

As in [V2]--[V4], the topology of the space $L(d) \setminus \Sigma(d)$ reduces
to that of the discriminant $\Sigma(d)$. The cohomology classes of this space
are deduced from a simplicial resolution of the discriminant. In this
resolution, over a map which glues together some $k_1$ points from the first
circle, $k_2$ points from the second and so on, exactly the simplex considered
in the definition of the complex $K (A; A_1, \dots , A_d)$ appears. This
simplex lies in the $k-1$st term of the natural filtration of the resolution,
while the union of its faces corresponding to the non-connected graphs lies in
the term $F_{k-2}.$

The simplest (of filtration 1) homotopy invariants obtained from this spectral
sequence are just the linking numbers of different components of the link; the
$p$ fold degrees and products of these linking numbers are the simplest
examples of the invariants of filtration $p$. \bigskip

\centerline{\bf References}
\medskip

[Al] V.I.~Arnold, {\it The cohomology ring of the group of colored braids},
Mat. Zametki, {\bf 5}(1969), 227-231; English translation: Math. Notes,
5(1969), 138-140.  \smallskip

[A2] V.I.~Arnold, {\it On some topological invariants of the algebraic
functions}, Trans. Moscow Math. Soc., {\bf 21}(1970), 27-46. \smallskip

[A3] V.I.~Arnold, {\it Spaces of functions with mild singularities}, Funct.
Anal. Appl., {\bf 23}:3 (1989), 169--177. \smallskip

[AGLV] V.I.~Arnold, V.V.~Gorjunov, O.V.~Ljashko, V.A.~Vassiliev, {\it
Singularities II: Classification and Applications}, Itogi nauki VINITI,
Fundamentalnyje napravlenija, {\bf 39}(1989) Moscow, VINITI; English
translation: Encycl. Math. Sci., {\bf 39}(1993), Berlin a.o.: Springer.
\smallskip

[BLY] A.~Bjorner, L.~Lovasz, A.~Yao, {\it Linear decision trees: volume
estimates and topological bounds}, Report No. 5 (1991/1992), Inst.
MittagLeffler (1991). \smallskip

[BW] A.~Bjorner, V.~Welker, {\it The homology of ``$k$-equal'' manifolds and
related partition lattices}. Preprint, 1992. \smallskip

[Br] E.~Brieskorn, {\it Sur les groupes de tresses (d'apres V.I. Arnold)}, Sem.
Bourbaki, 1971/72, No. 401, Springer Lect. Notes Math, (317) (1973), p. 21-44.
\smallskip

[D] P.~Deligne, {\it Theorie de Hodge}, II, III, Publ. Math. IHES, {\bf
40}(1970), 235 5-58, {\bf 44}(1972), 5-77. \smallskip

[FNRS] V.V.~Fock, N.A.~Nekrasov, A.A.~Rosly, K.G.~Selivanov, {\it What we think
about the higher dimensional Chern--Simons theories}, Preprint Inst. Theor. and
Experim. Phys., No. 70-91, (1991), Moscow. Published in Proceedings of First
international Sakharov conference on Physics, 1992. \smallskip

[F] J.~Folkman, {\it The homology group of a lattice.} J. Math. and Mech., {\bf
15}(1966), 631-636. \smallskip

[FV] D.B.~Fuchs, O.Ya.~Viro, {\it Introduction to the homotopy theory}, Itogi
nauki VINITI, Fundamentalnyje napravlenija, 24(1988), Moscow VINITI; English
translation in  {\it Topology II}, vol. 24 of Encycl. Math. Sci, 24, pages
1--93, Berlin a.o.: Springer, 2004.
\smallskip

[GM] M.~Goresky, R.~MacPherson, {\it Stratified Morse Theory}, Berlin a.o.:
Springer (1986). \smallskip

[GR] I.M.~Gelfand, G.L.~Rybnikov, {\it Algebraic and topological invariants of
oriented matroids}. Soviet Math. Doklady, {\bf 33}(1986), 573-577. \smallskip

[M] J.W.~Milnor, {\it Isotopy of links}. In: Algebraic Geometry and Topology,
Princeton, N.J., Princeton Univ. Press, 1957, 280-306. \smallskip

[N] N.A.~Nekrasov, {\it On the cohomology of the complement of the bifurcation
diagram of the singularity} $A_\mu$. Funct. Anal. and its Appl. {\bf 27}:4
(1993), 245--250. \smallskip

[O] P.~Orlik, {\it Introduction to Arrangements}, CBMS Lecture Notes, AMS {\bf
72}, 1989. \smallskip

[OS] P.~Orlik, L.~Solomon, {\it Combinatorics and topology of complements and
hyperplanes}, Invent. Math., {\bf 56}(1980), 167-189. \smallskip

[V1] V.A.~Vassiliev, {\it The topology of the complement of a plane
arrangement}, e-preprint, 1991. \smallskip

[V2] V.A.~Vassiliev, {\it Complements of Discriminants of Smooth Maps: Topology
and Applications}, AMS, Translations of Math. Monographs, ({\bf 98}) (1992),
Providence, R.I. \smallskip

[V3] V.A.~Vassiliev, {\it Homological invariants of knots: algorithms and
calculations}. Preprint Inst. Appplied Math. (90)(1990), Moscow (in Russian).
\smallskip

[V4] V.A.~Vassiliev, {\it Cohomology of Knot Spaces}. In: Theory of
Singularities and its Applications, V.I. Arnold, ed., AMS, Advances in Soviet
Math., {\bf 1}(1990), 23-69. \smallskip

[W] C.W.~Whitehead, {\it Recent Advances in Homotopy Theory}, Publ. AMS,
(1970).
\smallskip

[ZZ] G.M.~Ziegler, R.T.~Zivaljevich, {\it Homotopy type of arrangements via
diagrams of spaces}, Report No. {\bf 10}(1991/1992), Inst. Mittag-Leffler,
December 1991.
\end{document}